\documentclass[a4paper,10pt]{article}
\usepackage[T1]{fontenc}
\usepackage[utf8]{inputenc}
\usepackage{graphicx}
\usepackage{fancyhdr}
\usepackage{float}
\usepackage{xcolor}
\usepackage{authblk}
\usepackage{mathrsfs}
\usepackage{empheq}
\usepackage[hyphens]{url}
\usepackage{hyperref} 
\usepackage[]{breakurl}
\usepackage[]{amsthm}
 
\newcommand{\ninej}[9]{\left\{ \begin{array}{ccc}
            #1 & #2 & #3 \\
			#4 & #5 & #6 \\
			#7 & #8 & #9
			\end{array} \right\} }

\begin{document}

\huge

\begin{center}
Expression of special stretched $9j$ coefficients in terms of $_5F_4$ hypergeometric series
\end{center}

\vspace{0.5cm}

\large

\begin{center}
Jean-Christophe Pain$^{a,b,}$\footnote{jean-christophe.pain@cea.fr}
\end{center}

\normalsize

\begin{center}
\it $^a$CEA, DAM, DIF, F-91297 Arpajon, France\\
\it $^b$Universit\'e Paris-Saclay, CEA, Laboratoire Mati\`ere en Conditions Extr\^emes,\\
\it F-91680 Bruy\`eres-le-Ch\^atel, France\\
\end{center}

\vspace{0.5cm}

\begin{abstract}
The Clebsch-Gordan coefficients or Wigner $3j$ symbols are known to be proportional to a $_3F_2(1)$ hypergeometric series, and Racah $6j$ coefficients to a $_4F_3(1)$. In general, however, non-trivial $9j$ symbols can not be expressed as a $_5F_4$. In this letter, we show, using the Dougall-Ramanujan identity, that special stretched $9j$ symbols can be reformulated as $_5F_4(1)$ hypergeometric series.
\end{abstract}

\vspace{0.5cm}

\section{Introduction}\label{sec1}

It is well known that the Clebsch-Gordan coefficients or $3j$ symbols can be expressed by a $_3F_2(1)$ hypergeometric series. In the same way, $6j$ coefficients are proportional to a $_4F_3(1)$ \cite{Varshalovich1988}. The $9j$ angular-momentum coefficient plays a major role in atomic physics, for example because it characterizes the transformation from $LS$ to $jj$ couplings. The simplest known form for the $9j$ is due to Ali$\check{\mathrm{s}}$auskas, Jucys and Bandzaitis \cite{Alisauskas1972,Jucys1977,Rao1989}. It consists in a triple summation, which numerical implementation is faster and more accurate than the conventional single sum (over one product of three $6j$ coefficients) and double sum (over a product of $3j$ coefficients) \cite{Rosengren1998}. The Ali$\check{\mathrm{s}}$auskas-Jucys-Bandzaitis formula has been identified as a particular case of the triple hypergeometric functions of Lauricella-Saran-Srivastava \cite{Srivastava1967,Lauricella1893,Jeugt1994}. The $9j$ coefficient was shown not to belong, in general, to the $_pF_q$ family of hypergeometric functions \cite{Wu1973}. Wu investigated the class of hypergeometric functions of the Gel'fand type \cite{Wu1973,Gelfand1966}, being the Radon transforms of products of linear forms, but did not succeed. In general, the $9j$ symbol can therefore not be expressed as a $_5F_4(1)$. The question we try to answer in this letter is: ``are there particular $9j$ that can be expressed in terms of $_5F_4(1)$ (without summation)?''. Asked this way, the answer is obviously ``yes'', because if the arguments of the $9j$ are sufficiently degenerate, in the sense that one or more triads (columns or rows) are such that one of their argument is the sum of the two others, it will be a $6j$, or a $3j$, or even equal to 1, and we can then obviously write it as a $_5F_4(1)$ with particular arguments. In the present case, we add of course the additional constraint that the result (discarding the pre-factors that are products of factorials or factorial square roots) cannot be proportional to a $_4F_3(1)$, nor a $_3F_2(1)$. In short, the aim is to find a sufficiently non-trivial example.

A long time ago, Bandzaitis, Karosiene and Yutsis \cite{Bandzaitis1964} derived formulas for the stretched $9j$ coefficients. Sharp showed that there are five different doubly stretched coefficients and two triply stretched ones \cite{Sharp1967,Rao1994}. Some of them vanish, such as \cite{Zamick2011,Pain2011b}:
\begin{equation}
    \ninej{j}{j}{(2j-1)}{j}{j}{(2j-1)}{(2j-1)}{(2j-3)}{(4j-4)}=0.
\end{equation}
If one element is equal to unity \cite{Pain2011}, one gets a $6j$ or a linear combination of $6j$. If one argument is zero, the $9j$ reduces to a $6j$. If one triad is equal to $(1/2,1/2,1)$, the $9j$ reduces to a $3j$ \cite{Edmonds1957}. The following particular $9j$ with two degenerate triads
\begin{equation}
   \ninej{a}{b}{c}{d}{e}{f}{a+d}{b+e}{g}
\end{equation}
is also proportional to a $3j$ (see, e.g., Refs. \cite{Kleszyk2014} or \cite{Varshalovich1988}, Eq. (9) p. 354). One has also, for instance \cite{Rao1994}:
\begin{eqnarray}
    \ninej{a}{b}{c}{d}{e}{f}{a+d}{a+d+g}{g}&=&(-1)^{d-e+f}\frac{\eta(a+d+g,b,e)}{\eta(a,b,c)\eta(d,e,f)\eta(g,c,f)}\nonumber\\
    & &\times\left[\frac{(2a)!(2d)!(2g)!}{(2a+2d+1)(2a+2b+2g+1)!}\right]^{1/2},
\end{eqnarray}
where
\begin{equation}
\eta(a,b,c)=\left[\frac{(a-b+c)!(a+b-c)!(a+b+c+1)!}{(-a+b+c)!}\right]^{1/2}.
\end{equation}
In the case where one of the argument is zero, a $_4F_3$ can be obtained, for instance (see Eq. (14) in Ref. \cite{Rao1998}):
\clearpage
\begin{eqnarray}
    \ninej{a}{a}{0}{d}{e}{f}{g}{h}{f}&=&\frac{(-1)^{a+e+f+g}(-1)^{\beta_1+1}}{\left[(2a+1)(2f+1)\right]^{1/2}}\Delta(a,e,h)\Delta(f,g,h)\Delta(a,g,d)\Delta(f,e,d)\nonumber\\
    & &\times\frac{\Gamma(1-\beta_1)}{\left[\Gamma\left(1-\alpha_1,1-\alpha_2,1-\alpha_3,1-\alpha_4,\beta_2,\beta_3\right)\right]^{1/2}}\nonumber\\
    & &\times ~_4F_3\left[
\begin{array}{c}
\alpha_1,\alpha_2,\alpha_3,\alpha_4\\
\beta_1,\beta_2,\beta_3
\end{array};1\right],\nonumber\\
\end{eqnarray}
where $\alpha_1=h-a-e$, $\alpha_2=h-f-g$, $\alpha_3=d-a-g$, $\alpha_4=d-e-f$, $\beta_1=-a-e-f-g-1$, $\beta_2=d+h-e-g+1$ and $\beta_3=d+h-a-f+1$. The notation $\Gamma(x,y,z,...)$ is a shortcut for the product  $\Gamma(x)\Gamma(y)\Gamma(z)...$.

In section \ref{sec2}, we show, using the Dougall-Ramanujan identity, that special stretched $9j$ symbols can be expressed as a $_5F_4(1)$ hypergeometric series. 

\section{From the Dougall-Ramanujan relation to $_5F_4(1)$ representation of particular stretched $9j$ symbols}\label{sec2}

Varshalovich \emph{et al.} give (formula (13) p. 354) \cite{Varshalovich1988}:
\begin{eqnarray}\label{var}
    \ninej{a}{b}{a+b}{d}{e}{f}{e}{d}{a+b+f}&=&(-1)^{a+d-e}\frac{\Delta(a+b+f,e,d)}{\Delta(a,d,e)\Delta(b,e,d)\Delta(d,e,f)}\nonumber\\
    & &\times\left[\frac{(2a)!(2b)!(2f)!}{(2a+2b+1)(2a+2b+2f+1)!}\right]^{1/2}\nonumber\\
    & &\times\frac{(a+b+e+d+f+1)}{(a+e+d+1)(b+e+d+1)(d+e+f+1)}\nonumber\\
    & &\times \frac{(a+b+e+d+f)!(e-a+d)!(e-b+d)!(d+e-f)!}{(e+d-a-b-f)!(a+e+d)!(b+e+d)!(d+e+f)!},\nonumber\\
    & &
\end{eqnarray}
where
\begin{equation}
\Delta(a,b,c)=\left[\frac{(a+b-c)!(a-b+c)!(-a+b+c)!}{(a+b+c+1)!}\right]^{1/2}.
\end{equation}
Now one uses the following Dougall-Ramanujan identity for the well-poised $_5F_4$ series, which reads \cite{Dougall1907,Gasper2004}:
\begin{eqnarray}\label{doug}
~_5F_4\left[
\begin{array}{c}
\frac{n}{2}+1,n,-x,-y,-z\\
\frac{n}{2},x+n+1,y+n+1,z+n+1
\end{array};1\right]&=&\frac{\Gamma(x+n+1)\Gamma(y+n+1)}{\Gamma(n+1)\Gamma(x+y+z+n+1)}\nonumber\\
& &\times\frac{\Gamma(z+n+1)\Gamma(x+y+n+1)}{\Gamma(y+z+n+1)\Gamma(x+z+n+1)}.\nonumber\\
& &
\end{eqnarray}
Note that for $z=-n/2$, the latter expression reduces to the Dixon formula \cite{Dixon1902}. The right-hand side of the above equation can be re-expressed as the last term of the formula (\ref{var}), namely
\begin{equation}
\frac{(a+b+e+d+f)!(e-a+d)!(e-b+d)!(d+e-f)!}{(e+d-a-b-f)!(a+e+d)!(b+e+d)!(d+e+f)!}
\end{equation}
provided that one performs the change of variables $x=b+f$, $y=a+f$, $z=a+b$ and $n=e+d-a-b-f$. One thus gets
\begin{eqnarray}\label{newres}
    \ninej{a}{b}{a+b}{d}{e}{f}{e}{d}{a+b+f}&=&(-1)^{a+d-e}\frac{\Delta(a+b+f,e,d)}{\Delta(a,d,e)\Delta(b,e,d)\Delta(d,e,f)}\nonumber\\
    & &\times\left[\frac{(2a)!(2b)!(2f)!}{(2a+2b+1)(2a+2b+2f+1)!}\right]^{1/2}\nonumber\\
    & &\times\frac{(a+b+e+d+f+1)}{(a+e+d+1)(b+e+d+1)(d+e+f+1)}\nonumber\\
    & &\times~_5F_4\left[
\begin{array}{c}
\frac{e+d-a-b-f}{2}+1,e+d-a-b-f,-b-f,-a-f,-a-b\\
\frac{e+d-a-b-f}{2},e+d-a+1,e+d-b+1,e+d-f+1
\end{array};1\right],\nonumber\\
& &
\end{eqnarray}
which is the main result of the present work. As an example, one finds
\begin{equation}
\ninej{6}{10}{16}{14}{12}{8}{12}{14}{24}=\frac{13}{124062}\sqrt{\frac{1615}{7683753}},
\end{equation}
and the computation times are, using the \verb!AbsoluteTiming[]! function in Mathematica \cite{Mathematica}: 0.999675 s. for the Racah formula (a sum of product of three $6j$ symbols), 0.001472 s. for the formula (13), p. 354 section 10.8.3 of Ref. \cite{Varshalovich1988} and 0.000542 s. using the $_5F_4(1)$ hypergeometric-series expression (see Eq. (\ref{newres})). The new expression is thus the most efficient.

\section{Numerical implementation}

The interest of hypergeometric functions for the computation of $3nj$ coefficients was already pointed out. Wills \cite{Wills1971} arranged the series expansion of the Clebsch-Gordan coefficient into a nested form and suggested that a similar rearrangement was possible for the $6j$ coefficient, which was confirmed by Bretz \cite{Bretz1976}. We propose to use the same technique for the particular stretched $9j$ coefficients of Eq. (\ref{newres}). The basic idea is to compute the $_pF_q(1)$ using Horner's rule for polynomial evaluation, as
\begin{equation}\label{para}
~_pF_q\left[
\begin{array}{c}
\alpha_1,\cdots, \alpha_p,\\
\beta_1,\cdots,\beta_q
\end{array};z\right]=\left[1+\frac{\mathscr{A}_0}{\mathscr{B}_0}\left(z+\frac{\mathscr{A}_1}{\mathscr{B}_1}\left(z+\frac{\mathscr{A}_2}{\mathscr{B}_2}\left(\vphantom{\frac{\mathscr{A}_1}{\mathscr{B}_1}}z+\cdots\right)\right)\right)\right],
\end{equation} 
where
\begin{equation}
\mathscr{A}_i=\prod_{j=1}^p\left(\alpha_j+i\right)
\end{equation}
and
\begin{equation}
\mathscr{B}_i=(i+1)\prod_{k=1}^q\left(\beta_k+i\right).
\end{equation}
The Wills nested form avoids the numerical overflow issues due to the factorials and is rather fast \cite{Rao1981}. This is important in the case of large angular momenta, because the summations may contain large numbers of terms.

\section{Comments and precautions}

It is important, however, to keep in mind that $3j$ symbols play an essentially different role in representation theory than $6$ and $9j$ symbols. Moreover, $6j$ symbols are given by {\it balanced} (or {\it Saalsch\"utzian}) $_4F_3$ series, which means, using the notations of Eq. (\ref{para}), that 
\begin{equation}
\sum_{i=1}^4\beta_i=1+\sum_{i=1}^3\alpha_i.    
\end{equation}
Conversely, the sums appearing in the present paper are {\it well poised} $_5F_4$ series, which means that
\begin{equation}
    1+\alpha_1=\beta_1+\alpha_2=\beta_2+\alpha_3=\beta_3+\alpha_4.
\end{equation}
Although it makes sense to think of the $9j$ symbol as a generalization of the $6j$ symbol, the well-poised $_5F_4$ is not a generalization of the balanced $_4F_3$.

\section{Conclusion}

It seems that there is no way to express a general $9j$ (with 9 arbitrary parameters satisfying the required triangular inequalities), using a simple $_5F_4(1)$, without summation. In the present work, we found an example of doubly-stretched $9j$ symbol, which can be put in a form proportional to a $_5F_4(1)$ hypergeometric function. The example in question contains five free parameters out of nine, which is already important, but there may be another (less) particular $9j$ containing more free parameters, which may be expressed as a $_5F_4(1)$, times a relatively universal pre-factor.

\vspace{5mm}

{\bf Acknowledgements}

\vspace{5mm}

I would like to thank Robert Coquereaux, from Aix-Marseille University (France), for his help and useful comments in the course of this work. Valuable discussions with Joris Van der Jeugt, from Gent University (Belgium), are also acknowledged.

\section*{Appendix A: Remarks on the implementation of hypergeometric functions in a Computer Algebra System}\label{appA}

The evaluation, with a Computer Algebra System, of coefficients
\begin{equation}
	\ninej{a}{b}{a+b}{d}{e}{f}{e}{d}{a+b+f}
\end{equation}
when $e + d - a - b - f=0$, requires special attention. For instance, using Mathematica \cite{Mathematica}:
\vspace{5mm}
\begin{verbatim}
HypergeometricPFQ[{n/2+1, n, -x, -y, -z}, {n/2, x+n+1, y+n+1, z+n+1}, 1] 
/. {n -> 0, x -> 2, y -> 2, z -> 2} 
\end{verbatim}
\vspace{5mm}
gives the value 1, but
\vspace{5mm}

\begin{verbatim} 
HypergeometricPFQ[{n/2+1, n, -x, -y, -z}, {n/2, x+n+1, y+n+1, z+n+1}, 1] 
/. {x -> 2, y -> 2, z -> 2} /. {n -> 0}
\end{verbatim}
\vspace{5mm}
which should provide the same result, gives the value 5/12. The latter value is actually the correct result (and the former is by the way incompatible with the Dougall-Ramanujan identity (\ref{doug})).

%\section*{Declarations}

%\begin{itemize}
%\item No funding was received for conducting this study.
%\item There is no conflict of interest.
%\item Ethics approval and consent to participate
%\item Consent for publication
%\item No data was used for this research. 
%\item Materials availability
%\item Code availability 
%\item Author contribution
%\end{itemize}

\end{document}